\documentclass[11pt]{amsart}
\usepackage{latexsym}
\usepackage{amsmath}
\usepackage{amssymb}

\newcommand{\eproof}{\mbox{\ }\hfill $\Box$ \par \vskip 10pt}

\newtheorem{Theorem}{Theorem}[section]
\newtheorem{lemma}[Theorem]{Lemma}
\newtheorem{prop}[Theorem]{Proposition}

\numberwithin{equation}{section}

\def\cal{\mathcal}

\begin{document}

\title[Improved resolvent bounds]{Improved resolvent bounds for radial potentials}

\author[G. Vodev]{Georgi Vodev}

\address {Universit\'e de Nantes, Laboratoire de Math\'ematiques Jean Leray, 2 rue de la Houssini\`ere, BP 92208, 44322 Nantes Cedex 03, France}
\email{Georgi.Vodev@univ-nantes.fr}

\date{}

\begin{abstract} We prove semiclassical resolvent estimates for the
Schr\"odinger operator in $\mathbb{R}^d$, $d\ge 3$, with real-valued radial potentials $V\in L^\infty(\mathbb{R}^d)$.
In particular, we show that if $V(x)={\cal O}\left(\langle x\rangle^{-\delta}\right)$ with $\delta>2$, then the resolvent bound
is of the form $\exp\left(Ch^{-4/3}\right)$ with some constant $C>0$. We also get resolvent bounds when $1<\delta\le 2$. 
For slowly decaying $\alpha$ - H\"older potentials we get better resolvent bounds of the form $\exp\left(Ch^{-4/(\alpha+3)}\right)$.

\quad

\end{abstract} 

\maketitle

\setcounter{section}{0}
\section{Introduction and statement of results}

The aim of this work is to improve the recent results in \cite{kn:GS}, \cite{kn:V2}, \cite{kn:V4} concerning
the semiclassical behavior of the resolvent of the Schr\"odinger operator
$$P(h)=-h^2\Delta+V(x)$$
where $0<h\ll 1$ is a semiclassical parameter, $\Delta$ is the negative Laplacian in 
$\mathbb{R}^d$, $d\ge 3$, and $V\in L^\infty(\mathbb{R}^d)$ is a real-valued short-range potential satisfying the condition
\begin{equation}\label{eq:1.1}
|V(x)|\le C(|x|+1)^{-\delta}
\end{equation}
where $C>0$ and $\delta>1$ are some constants.
 More precisely, we are interested in bounding  
the quantity
$$g_s^\pm(h,\varepsilon):=\log\left\|(|x|+1)^{-s}(P(h)-E\pm i\varepsilon)^{-1}(|x|+1)^{-s}
\right\|_{L^2(\mathbb{R}^d)\to L^2(\mathbb{R}^d)}$$
from above by an explicit function of $h$, independent of $\varepsilon$. 
Here $0<\varepsilon<1$, $s>1/2$ is independent of $h$ and $E>0$ is a fixed energy level independent of $h$.
When $\delta>2$ it has been proved in \cite{kn:GS} that
\begin{equation}\label{eq:1.2}
g_s^\pm(h,\varepsilon)\le Ch^{-4/3}\log(h^{-1}).
\end{equation}
The bound (\ref{eq:1.2}) was previously proved in \cite{kn:KV} and \cite{kn:S3} when $d\ge 2$ for compactly supported potentials,
and in \cite{kn:V1} when $\delta>3$ and $d\ge 3$. It was also shown in \cite{kn:V3} that (\ref{eq:1.2}) still holds for
more general asymptotically Euclidean manifolds. In the present paper we show that, if $d\ge 3$, the logarithmic term in the righ-hand side of
(\ref{eq:1.2}) can be removed for potentials $V$ depending only on the radial variable $r=|x|$. We also improve significantly the bound
\begin{equation}\label{eq:1.3}
g_s^\pm(h,\varepsilon)\le 
 Ch^{-\frac{2\delta+5}{3(\delta-1)}}\left(\log(h^{-1})\right)^{\frac{1}{\delta-1}}.
\end{equation}
proved in \cite{kn:V2} for $1<\delta\le 3$ and $d\ge 3$. More precisely, we have the following

\begin{Theorem} Let $d\ge 3$ and suppose that the potential $V$ depends only on the radial variable. If $V$ 
satisfies (\ref{eq:1.1}) with $\delta>2$,  
then there exist constants $C>0$ and $h_0>0$ independent of $h$ and $\varepsilon$ but depending on $s$, $E$, such 
that the bound
\begin{equation}\label{eq:1.4}
g_s^\pm(h,\varepsilon)\le Ch^{-4/3}
\end{equation}
holds for all $0<h\le h_0$. If $V$ 
satisfies (\ref{eq:1.1}) with $1<\delta\le 2$,  
then we have the bound
\begin{equation}\label{eq:1.5}
g_s^\pm(h,\varepsilon)\le 
 Ch^{-\frac{2\delta}{2\delta-1}}\left(\log(h^{-1})\right)^{\frac{\delta+1}{2\delta-1}}.
\end{equation}
If the potential satisfies the condition 
\begin{equation}\label{eq:1.6}
|V(r)|\le \widetilde C(r+1)^{-1}(\log(r+2))^{-\rho}
\end{equation}
with some constants $\widetilde C>0$ and $\rho>1$, 
then we have the bound
\begin{equation}\label{eq:1.7}
g_s^\pm(h,\varepsilon)\le Ch^{-2}.
\end{equation}
\end{Theorem}

Note that when $d=1$ we have much better resolvent bounds. Indeed, it has been proved in 
 \cite{kn:DS} that for $V\in L^1(\mathbb{R})$ we have 
the bound
\begin{equation}\label{eq:1.8}
g_s^\pm(h,\varepsilon)\le Ch^{-1}.
\end{equation}
The bound (\ref{eq:1.8}) is proved in \cite{kn:V4} (see also \cite{kn:D}) when $d\ge 3$ for slowly decaying Lipschitz potentials $V$ 
with respect to the radial variable $r$ and satisfying the conditions
\begin{equation}\label{eq:1.9}
V(x)\le p(|x|)
\end{equation}
where $p(r)>0$, $r\ge 0$, is a decreasing function such that 
$p(r)\to 0$ as $r\to\infty$, and 
\begin{equation}\label{eq:1.10}
\partial_rV(x)\le C_1(|x|+1)^{-\delta_1}
\end{equation}
where $C_1>0$ and $\delta_1>1$ are some constants. When $d=2$ the bound (\ref{eq:1.8}) is proved in \cite{kn:S2}
for potentials which are Lipschitz with respect to the space variable $x$. Under this condition, 
the bound (\ref{eq:1.8}) is extended in \cite{kn:V4} to general exterior domains and all dimensions $d\ge 2$. 
Note that the bound (\ref{eq:1.8}) was first proved for smooth potentials
in \cite{kn:B2}. 
A high-frequency analog of (\ref{eq:1.8}) on Riemannian manifolds was proved in 
\cite{kn:B1} and \cite{kn:CV}. It was also showed in \cite{kn:DDZ} that the bound (\ref{eq:1.8}) is optimal for smooth potentials. 

The bound (\ref{eq:1.2}) has been recently improved in \cite{kn:GS}, \cite{kn:V4} for H\"older potentials $V\in C_\beta^\alpha(\overline{\mathbb{R}^+})$ with respect to the radial variable and satisfying the condition
(\ref{eq:1.9}) as well. Hereafter, given $0<\alpha<1$ and $\beta>0$, 
the space $C_\beta^\alpha(\overline{\mathbb{R}^+})$ denotes the set of all H\"older functions $a$ such that
 $$\sup_{r'\ge 0:\,0<|r-r'|\le 1}\frac{|a(r)-a(r')|}{|r-r'|^\alpha}\le C(r+1)^{-\beta},\quad\forall r\in \overline{\mathbb{R}^+},$$
 for some constant $C>0$. Indeed, it is shown in \cite{kn:GS} (with $\beta>3$) and in \cite{kn:V4} (with $\beta=4$) 
 that in this case we have the bound
 \begin{equation}\label{eq:1.11}
g_s^\pm(h,\varepsilon)\le Ch^{-4/(\alpha+3)}\log(h^{-1}).
\end{equation}
In \cite{kn:V4} the bound (\ref{eq:1.11}) is also extended to exterior domains and all dimensions $d\ge 2$ for potentials 
which are $\alpha$ - H\"older with respect to the space variable. In the present paper we show that the bound (\ref{eq:1.11})
can be improved for radial potentials. We have the following

\begin{Theorem} Let $d\ge 3$ and suppose that the potential $V$ depends only on the radial variable and satisfies (\ref{eq:1.9}). If 
$V\in C_3^\alpha(\overline{\mathbb{R}^+})$,  
then there exist constants $C>0$ and $h_0>0$ independent of $h$ and $\varepsilon$ but depending on $s$, $E$ and the function $p$, such 
that the bound
\begin{equation}\label{eq:1.12}
g_s^\pm(h,\varepsilon)\le Ch^{-4/(\alpha+3)}
\end{equation}
holds for all $0<h\le h_0$. If $V\in C_\beta^\alpha(\overline{\mathbb{R}^+})$ with $2<\beta<3$, then we have the bound
\begin{equation}\label{eq:1.13}
g_s^\pm(h,\varepsilon)\le 
 Ch^{-k}\left(\log(h^{-1})\right)^q,
\end{equation}
where 
$$k=\frac{2\alpha\beta-6\alpha+4}{2\alpha\beta-5\alpha+3},\quad q=\frac{\alpha(3-\beta)}{2\alpha\beta-5\alpha+3}.$$
If $V\in C_\beta^\alpha(\overline{\mathbb{R}^+})$ with $1<\beta\le 2$, then the bound (\ref{eq:1.13}) holds with 
 $$k=\frac{2\beta-2\alpha}{2\beta-\alpha-1},\quad q=\frac{\beta+1}{2\beta-\alpha-1}.$$
\end{Theorem}

\noindent
{\bf Remark 1.} It is easy to see from the proof that the resolvent bounds in the above theorems still hold if we add to the potential $V$
a radial real-valued Lipschitz function satisfying (\ref{eq:1.9}) and (\ref{eq:1.10}).

\noindent
{\bf Remark 2.} The conclusions of the above theorems remain valid for the Dirichlet self-adjoint realisation of the operator
$P(h)$ on the Hilbert space $L^2(\mathbb{R}^d\setminus B)$, where $B=\{x\in \mathbb{R}^d:|x|\le r_0\}$, $r_0>0$. 
In fact, in this case the proof also works when $d=2$. The reason why the proof does not work on the whole $\mathbb{R}^2$
 is that in this case the effective potential is negative with a strong singularity at $r=0$. 

Resolvent bounds like those in the above theorems are usually proved by means of Carleman estimates with suitable
phase functions. Since the resolvent bound is determined by the magnitude of the phase, the aim is to find as small phase
as possible in order to get as good bound as possible. In the case of general potentials this proves to be a quite
delicate problem which is not easy to handle with. Our proof in the present paper is very different from
the proofs in the papers mentioned above. We take advantage 
of the assumption that the potential is radial in order to reduce our
$d$-dimensional resolvent estimate to infinitely many one dimensional resolvent estimates depending on an additional
parameter, which is the eigenvalues of the Laplace-Beltrami operator on the unit sphere (see Section 2). 
This allows us to use simpler Carleman estimates with better phases, and therefore we can get better bounds.
We expect that the bound (\ref{eq:1.4}) is optimal for $L^\infty$ potentials, but to our best knoweledge
showing this remains an open problem. We also expect that (\ref{eq:1.4}) holds
 for general potentials satisfying (\ref{eq:1.1}) with $\delta>2$, without assuming that $V$
 is radial. This also is an open problem in the context of Euclidean spaces. 
Note that (\ref{eq:1.4}) is proved in \cite{kn:V3} for a quite large class of
asymptotically hyperbolic manifolds with $L^\infty$ potentials decaying sufficiently fast.
Finally, the bounds (\ref{eq:1.4}) and (\ref{eq:1.12}) imply a better local energy decay
for the wave equation with a radial refraction index than that one in \cite{kn:S1}
and \cite{kn:V4}.

\section{Preliminaries}

To get our resolvent bounds we will use the following

\begin{lemma} Let $s>1/2$ and suppose that for all functions
$f\in H^2(\mathbb{R}^d)$ such that 
$$(|x|+1)^{s}(P(h)-E\pm i\varepsilon)f\in L^2(\mathbb{R}^d)$$
 we have the estimate 
 \begin{equation}\label{eq:2.1}
\|(|x|+1)^{-s}f\|^2_{L^2(\mathbb{R}^d)}
\le M\|(|x|+1)^{s}(P(h)-E\pm i\varepsilon)f\|^2_{L^2(\mathbb{R}^d)}+M\varepsilon\|f\|^2_{L^2(\mathbb{R}^d)}
\end{equation}
with some $M>0$ independent of $\varepsilon$ and $f$. Then we have the resolvent bound
\begin{equation}\label{eq:2.2}
g_s^\pm(h,\varepsilon)\le \log(M+1).
\end{equation}
\end{lemma}

{\it Proof.}  Since the operator $P(h)$ is symmetric, we have
$$
\varepsilon\|f\|^2_{L^2}=\pm{\rm Im}\,\langle (P(h)-E\pm i\varepsilon)f,f\rangle_{L^2}$$
$$\le (2M)^{-1}\|(|x|+1)^{-s}f\|^2_{L^2}+\frac{M}{2}\|(|x|+1)^{s}(P(h)-E\pm i\varepsilon)f\|^2_{L^2}$$
which can be rewritten in the form
\begin{equation}\label{eq:2.3}
M\varepsilon\|f\|^2_{L^2}\le \frac{1}{2}\|(|x|+1)^{-s}f\|^2_{L^2}+
\frac{M^2}{2}\|(|x|+1)^{s}(P(h)-E\pm i\varepsilon)f\|^2_{L^2}.
\end{equation}
By (\ref{eq:2.1}) and (\ref{eq:2.3}) we get
\begin{equation}\label{eq:2.4}
\|(|x|+1)^{-s}f\|_{L^2}\le (M+1)\|(|x|+1)^{s}(P(h)-E\pm i\varepsilon)f\|_{L^2}
\end{equation}
which clearly implies (\ref{eq:2.2}).
\eproof

We will now use that the potential is radial to reduce the estimate (\ref{eq:2.1}) to infinitely many similar estimates on $\mathbb{R}^+$.
To this end we will write the operator $P(h)$ in 
polar coordinates $(r,w)\in\mathbb{R}^+\times\mathbb{S}^{d-1}$, $r=|x|$, $w=x/|x|$ and we will use that $L^2(\mathbb{R}^d)=L^2(\mathbb{R}^+\times\mathbb{S}^{d-1}, r^{d-1}drdw)$. We have the identity
\begin{equation}\label{eq:2.5}
 r^{(d-1)/2}\Delta  r^{-(d-1)/2}=\partial_r^2+\frac{\widetilde\Delta_w}{r^2}
\end{equation}
where $\widetilde\Delta_w=\Delta_w-\frac{1}{4}(d-1)(d-3)$ and $\Delta_w$ denotes the negative Laplace-Beltrami operator
on $\mathbb{S}^{d-1}$. Set $v=r^{(d-1)/2}f$ and
$${\cal P}^\pm(h)=r^{(d-1)/2}(P(h)-E\pm i\varepsilon)r^{-(d-1)/2}.$$
Using (\ref{eq:2.5}) we can write the operator ${\cal P}^\pm(h)$ in the coordinates $(r,w)$ as follows
$${\cal P}^\pm(h)={\cal D}_r^2+\frac{\Lambda_w}{r^2}+V(r)-E\pm i\varepsilon $$
where we have put ${\cal D}_r=-ih\partial_r$ and $\Lambda_w=-h^2\widetilde\Delta_w$. Clearly, the estimate (\ref{eq:2.1}) can be rewritten in the form
\begin{equation}\label{eq:2.6}
\|(r+1)^{-s}v\|^2_{L^2(\mathbb{R}^+\times\mathbb{S}^{d-1})}
\le M\|(r+1)^{s}{\cal P}^\pm(h)v\|^2_{L^2(\mathbb{R}^+\times\mathbb{S}^{d-1})}+M\varepsilon\|v\|^2_{L^2(\mathbb{R}^+\times\mathbb{S}^{d-1})}.
\end{equation}
Let $\lambda_j\ge 0$ be the eigenvalues of $-\Delta_w$ repeated with the multiplicities and let $e_j\in L^2(\mathbb{S}^{d-1})$
be the corresponding eigenfunctions. Set
$$\nu=h\sqrt{\lambda_j+\frac{1}{4}(d-1)(d-3)}$$
and
$$v_j(r)=\langle v(r,\cdot),e_j\rangle_{L^2(\mathbb{S}^{d-1})},$$ 
$$Q^\pm_\nu(h)={\cal D}_r^2+\frac{\nu^2}{r^2}+V(r)-E\pm i\varepsilon.$$
Thus we can write 
$$v=\sum_j v_je_j,$$
$${\cal P}^\pm(h)v=\sum_jQ^\pm_\nu(h)v_j e_j,$$
so we have the identities
$$\|v\|^2_{L^2(\mathbb{R}^+\times\mathbb{S}^{d-1})}=\sum_j\|v_j\|^2_{L^2(\mathbb{R}^+)},$$
$$\|(r+1)^{-s}v\|^2_{L^2(\mathbb{R}^+\times\mathbb{S}^{d-1})}=\sum_j\|(r+1)^{-s}v_j\|^2_{L^2(\mathbb{R}^+)},$$
$$\|(r+1)^{s}{\cal P}^\pm(h)v\|^2_{L^2(\mathbb{R}^+\times\mathbb{S}^{d-1})}=\sum_j\|(r+1)^{s}Q_\nu^\pm(h)v_j\|^2_{L^2(\mathbb{R}^+)}.$$
We have the following

\begin{lemma} Let $s>1/2$ and suppose that for all $\nu$ the estimates
\begin{equation}\label{eq:2.7}
\|(r+1)^{-s}u\|^2_{L^2(\mathbb{R}^+)}
\le M_\nu\|(r+1)^{s}Q_\nu^\pm(h)u\|^2_{L^2(\mathbb{R}^+)}$$
$$+M_\nu\varepsilon\|u\|^2_{L^2(\mathbb{R}^+)}
+M_\nu\varepsilon\|{\cal D}_ru\|^2_{L^2(\mathbb{R}^+)}
\end{equation}
hold for every $u\in H^2(\mathbb{R}^+)$
such that $u(0)=0$ and $(r+1)^{s}Q_\nu^\pm(h)u\in L^2(\mathbb{R}^+)$, with $M_\nu>0$ independent of $\varepsilon$ and $u$.  
Then the estimate (\ref{eq:2.6}) holds with 
$$M=\left(2+E+\|V\|_{L^\infty}\right)\max_{\nu^2\in{\rm spec}\,\Lambda_w}M_\nu.$$
\end{lemma}

{\it Proof.} We integrate by parts to obtain
$${\rm Re}\int_0^\infty {\cal D}_r^2u\overline{u}dr=\int_0^\infty |{\cal D}_ru|^2dr,$$
which leads to 
$${\rm Re}\int_0^\infty Q_\nu^\pm(h)u\overline{u}dr=\int_0^\infty |{\cal D}_ru|^2dr+\nu^2\int_0^\infty r^{-2}|u|^2dr
+\int_0^\infty (V(r)-E)|u|^2dr$$
$$\ge \int_0^\infty |{\cal D}_ru|^2dr
-\left(E+\|V\|_{L^\infty}\right)\int_0^\infty |u|^2dr.$$
This implies 
\begin{equation}\label{eq:2.8}
\int_0^\infty |{\cal D}_ru|^2dr\le \int_0^\infty |Q_\nu^\pm(h)u|^2dr+\left(1+E+\|V\|_{L^\infty}\right)\int_0^\infty |u|^2dr.
\end{equation}
Combining (\ref{eq:2.7}) and (\ref{eq:2.8}) we get
\begin{equation}\label{eq:2.9}
\|(r+1)^{-s}u\|^2_{L^2(\mathbb{R}^+)}
\le 2M_\nu\|(r+1)^{s}Q_\nu^\pm(h)u\|^2_{L^2(\mathbb{R}^+)}$$
$$+\left(2+E+\|V\|_{L^\infty}\right)M_\nu\varepsilon\|u\|^2_{L^2(\mathbb{R}^+)}.
\end{equation}
Applying (\ref{eq:2.9}) with $u=v_j$ and summing up all the inequalities clearly lead to (\ref{eq:2.6}) with the desired value of $M$.
\eproof

Thus we reduce our problem to proving estimates like (\ref{eq:2.7}) with as good bounds $M_\nu$ as possible. 
This will be carried out in the next sections.  

\section{Bounding $M$ for $L^\infty$ potentials}

We will first prove the following

\begin{prop} Let $V\in L^1(\mathbb{R}^+)\cap L^\infty(\mathbb{R}^+)$. Then the estimate (\ref{eq:2.7}) holds for all 
$\nu$ with $M_\nu=e^{C(\nu+1)/h}$, where $C>0$ is a constant independent of $\nu$ and $h$.
\end{prop}

{\it Proof.} We will first consider the simplier case when $\nu=0$. Note that this may happen only when $d=3$. Set
$$F(r)=E|u(r)|^2+|{\cal D}_ru(r)|^2$$
and observe that the first derivative of $F$ satisfies the identity
$$F'(r)=2h^{-1}{\rm Im}\,Vu\overline{{\cal D}_ru}-2h^{-1}{\rm Im}\,Q_0^\pm(h)u\overline{{\cal D}_ru}\pm2\varepsilon h^{-1}{\rm Re}\,u\overline{{\cal D}_ru}.$$
We have  
$$-F'(r)\le h^{-1}|V|\left(|u|^2+|{\cal D}_ru|^2\right)+h^{-1}\gamma(r+1)^{-2s}|{\cal D}_ru|^2$$
$$+h^{-1}\gamma^{-1}(r+1)^{2s}|Q_0^\pm(h)u|^2
+\varepsilon h^{-1}\left(|u|^2+|{\cal D}_ru|^2\right)$$
for any $\gamma>0$. Let $s>1/2$ and $\eta=\frac{1}{2}\min\{1,E\}$. Set $\mu=e^{\psi/h}$, where
$$\psi(r)=\int_0^r\left(\eta^{-1}|V(\sigma)|+(\sigma+1)^{-2s}\right)d\sigma\le \eta^{-1}\|V\|_{L^1}+(2s-1)^{-1}.$$
We have $|V(r)|\le \eta\psi'(r)$. Using this together with the identity $\mu=h\mu'/\psi'$, we obtain
$$-(\mu F)'=-\mu'F-\mu F'\le -\mu'F+
\eta\mu'\left(|u|^2+|{\cal D}_ru|^2\right)+h^{-1}\gamma\mu(r+1)^{-2s}|{\cal D}_ru|^2$$
$$+h^{-1}\gamma^{-1}\mu(r+1)^{2s}|Q_0^\pm(h)u|^2
+\varepsilon h^{-1}\mu\left(|u|^2+|{\cal D}_ru|^2\right)$$
$$\le -\frac{1}{2}\mu'F+h^{-1}\gamma\mu(r+1)^{-2s}|{\cal D}_ru|^2$$
$$+h^{-1}\gamma^{-1}\mu(r+1)^{2s}|Q_0^\pm(h)u|^2
+\varepsilon h^{-1}\mu\left(|u|^2+|{\cal D}_ru|^2\right).$$
Integrating this inequality gives
$$0\le F(0)=-\int_0^\infty(\mu F)'\le -\frac{1}{2}\int_0^\infty\mu'F+h^{-1}\gamma\mu_0\int_0^\infty(r+1)^{-2s}|{\cal D}_ru|^2$$
$$+h^{-1}\gamma^{-1}\mu_0\int_0^\infty(r+1)^{2s}|Q_0^\pm(h)u|^2
+\varepsilon h^{-1}\mu_0\int_0^\infty\left(|u|^2+|{\cal D}_ru|^2\right),$$
where $\mu_0=\max\mu\le e^{C/h}$, $C>0$. 
Using that $\mu'\ge\psi'/h\ge h^{-1}(r+1)^{-2s}$ we deduce from the above inequality
$$\frac{1}{2}\int_0^\infty(r+1)^{-2s}F\le\frac{h}{2}\int_0^\infty\mu'F\le \gamma\mu_0\int_0^\infty(r+1)^{-2s}|{\cal D}_ru|^2$$
$$+\gamma^{-1}\mu_0\int_0^\infty(r+1)^{2s}|Q_0^\pm(h)u|^2
+\varepsilon \mu_0\int_0^\infty\left(|u|^2+|{\cal D}_ru|^2\right).$$
Now we take $\gamma=(3\mu_0)^{-1}$ so that we can absorbe the first term in the right-hand side of the above inequality. Thus we get
the estimate 
$$\frac{1}{6}\int_0^\infty(r+1)^{-2s}F\le
3\mu_0^2\int_0^\infty(r+1)^{2s}|Q_0^\pm(h)u|^2
+\varepsilon\mu_0\int_0^\infty\left(|u|^2+|{\cal D}_ru|^2\right),$$
which clearly implies the desired bound.

Consider now the case $\nu>0$. Then $\nu\ge h\nu_0$ with some constant $\nu_0>0$. Let $\phi_j\in C^\infty(\mathbb{R})$, $j=0,1,2,$ be 
real-valued functions such that $0\le\phi_j\le 1$, $\phi'_j\ge 0$, $\phi_0(\sigma)=0$ for $\sigma\le 1/3$, 
$\phi_0(\sigma)=1$ for $\sigma\ge 1/2$, $\phi_1(\sigma)=0$ for $\sigma\le 1$, 
$\phi_1(\sigma)=1$ for $\sigma\ge 2$, $\phi_2(\sigma)=0$ for $\sigma\le 3$, 
$\phi_2(\sigma)=1$ for $\sigma\ge 4$. Set $\kappa=4\sqrt{1+E+\|V\|_{L^\infty}}$ and $V_\nu=\nu^2\phi_0(r/\kappa\nu)r^{-2}+V$,
$u_1=\phi_1(r/\kappa\nu)u$, $u_2=(1-\phi_2)(r/\kappa\nu)u$. Observe that 
$$(\nu^2r^{-2}+V)u_1=V_\nu u_1$$
and $V_\nu\in L^1$ with norm $\|V_\nu\|_{L^1}={\cal O}(\nu+1)$. As above, we are going to bound from below the first derivative
of the function
$$F_1(r)=E|u_1(r)|^2+|{\cal D}_ru_1(r)|^2.$$
We have the identity
$$F'_1(r)=2h^{-1}{\rm Im}\,V_\nu u_1\overline{{\cal D}_ru_1}-2h^{-1}{\rm Im}\,Q_\nu^\pm(h)u_1\overline{{\cal D}_ru_1}\pm2\varepsilon h^{-1}{\rm Re}\,u_1\overline{{\cal D}_ru_1}$$
$$=2h^{-1}{\rm Im}\,V_\nu u_1\overline{{\cal D}_ru_1}-2h^{-1}{\rm Im}\,\phi_1Q_\nu^\pm(h)u\overline{{\cal D}_ru_1}\pm2\varepsilon h^{-1}{\rm Re}\,u_1\overline{{\cal D}_ru_1}+T(r),$$
where
$$T(r)=-2h^{-1}{\rm Im}\,[{\cal D}_r^2,\phi_1(r/\kappa\nu)]u\overline{{\cal D}_r\phi_1(r/\kappa\nu)u}$$
$$=2(\kappa\nu)^{-1}\phi_1\phi'_1|{\cal D}_ru|^2+h(\kappa\nu)^{-2}(\phi_1\phi''_1+2\phi'^2_1){\rm Im}\,u\overline{{\cal D}_ru}
+h^2(\kappa\nu)^{-3}\phi'_1\phi''_1|u|^2$$
$$\ge -{\cal O}\left(h\nu^{-2}\right)(\phi'_1+|\phi''_1|)\left(|u|^2+|{\cal D}_ru|^2\right).$$
Thus we obtain
$$-F'_1(r)\le h^{-1}|V_\nu|\left(|u_1|^2+|{\cal D}_ru_1|^2\right)+h^{-1}\gamma(r+1)^{-2s}|{\cal D}_ru_1|^2$$
$$+h^{-1}\gamma^{-1}(r+1)^{2s}|Q_\nu^\pm(h)u|^2
+\varepsilon{\cal O}\left(h^{-1}\right)\left(|u|^2+|{\cal D}_ru|^2\right)$$
$$+{\cal O}\left(h\nu^{-2}\right)(\phi'_1+|\phi''_1|)\left(|u|^2+|{\cal D}_ru|^2\right)$$
for any $\gamma>0$. Let $1/2<s\le 1$, let $\eta$ be as above and let $\lambda\gg 1$ be a large parameter independent
of $h$ and $\nu$ to be fixed later on. Set $\mu=e^{\psi/h}$ with
$$\psi(r)=\int_0^r\left(\eta^{-1}|V_\nu(\sigma)|+\lambda(\sigma+1)^{-2s}\right)d\sigma\le \eta^{-1}\|V_\nu\|_{L^1}+\lambda(2s-1)^{-1}
\lesssim \nu+1.$$
Integrating by parts and absorbing the term involving the potential in the same way as above we arrive at the inequality
$$\frac{1}{2}\int_0^\infty\psi'\mu F_1\lesssim\gamma\mu_\nu\int_0^\infty(r+1)^{-2s}|{\cal D}_ru_1|^2$$
$$+\gamma^{-1}\mu_\nu\int_0^\infty(r+1)^{2s}|Q_\nu^\pm(h)u|^2$$
$$+\varepsilon\mu_\nu\int_0^\infty\left(|u|^2+|{\cal D}_ru|^2\right)$$
$$+h^2\nu^{-2}\mu(2\kappa\nu)\int_{\kappa\nu}^{2\kappa\nu}\left(|u|^2+|{\cal D}_ru|^2\right),$$
where $\mu_\nu=\max\mu\le e^{C(\nu+1)/h}$, $C>0$. 
On the othe hand, we have 
$$\int_0^\infty\psi'\mu F_1\ge \lambda\mu(3\kappa\nu)\int_{3\kappa\nu}^{4\kappa\nu}(r+1)^{-2s}F_1\ge 
\lambda(4\kappa\nu+1)^{-2}\mu(3\kappa\nu)\int_{3\kappa\nu}^{4\kappa\nu}F_1.$$
Combining both inequalities and using that $\mu(3\kappa\nu)>\mu(2\kappa\nu)>1$ yield
$$\int_{3\kappa\nu}^{4\kappa\nu}F_1\lesssim\gamma(\nu+1)^2\mu_\nu\int_0^\infty(r+1)^{-2s}|{\cal D}_ru_1|^2$$
$$+\gamma^{-1}(\nu+1)^2\mu_\nu\int_0^\infty(r+1)^{2s}|Q_\nu^\pm(h)u|^2$$
$$+\varepsilon(\nu+1)^2\mu_\nu\int_0^\infty\left(|u|^2+|{\cal D}_ru|^2\right)$$
$$+\lambda^{-1}h^2\nu^{-2}(\nu+1)^2
\int_{\kappa\nu}^{2\kappa\nu}\left(|u|^2+|{\cal D}_ru|^2\right).$$
Observe now that $h^2\nu^{-2}(\nu+1)^2\lesssim 1$ as long as $\nu\ge h\nu_0$. Therefore, we can rewrite the above inequality in the form
\begin{equation}\label{eq:3.1}
\int_{3\kappa\nu}^{4\kappa\nu}\left(|u|^2+|{\cal D}_ru|^2\right)\lesssim\gamma(\nu+1)^2\mu_\nu\int_0^\infty(r+1)^{-2s}|{\cal D}_ru_1|^2$$
$$+\gamma^{-1}(\nu+1)^2\mu_\nu\int_0^\infty(r+1)^{2s}|Q_\nu^\pm(h)u|^2$$
$$+\varepsilon(\nu+1)^2\mu_\nu\int_0^\infty\left(|u|^2+|{\cal D}_ru|^2\right)$$
$$+\lambda^{-1}\int_{\kappa\nu}^{2\kappa\nu}\left(|u|^2+|{\cal D}_ru|^2\right).
\end{equation}
On the other hand, the choice of $\kappa$
guarantees the inequality
$$(\nu^2r^{-2}+V(r)-E)|u_2|^2\ge |u_2|^2.$$
Therefore, integrating by parts we obtain
$${\rm Re}\int_0^\infty Q_\nu^\pm(h)u_2\overline{u_2}=\int_0^\infty|{\cal D}_ru_2|^2+\int_0^\infty(\nu^2r^{-2}+V(r)-E)|u_2|^2$$
$$\ge \int_0^\infty|{\cal D}_ru_2|^2+\int_0^\infty |u_2|^2,$$
which implies
$$\int_0^\infty\left(|u_2|^2+|{\cal D}_ru_2|^2\right)\le \int_0^\infty|Q_\nu^\pm(h)u_2|^2$$
$$\le\int_0^\infty|Q_\nu^\pm(h)u|^2+\int_0^\infty\left|[{\cal D}_r^2,\phi_2(r/\kappa\nu)]u\right|^2$$
$$\lesssim\int_0^\infty|Q_\nu^\pm(h)u|^2+\int_{3\kappa\nu}^{4\kappa\nu}\left(|u|^2+|{\cal D}_ru|^2\right).$$
In particular, this inequality yields
\begin{equation}\label{eq:3.2}
\int_{\kappa\nu}^{2\kappa\nu}\left(|u|^2+|{\cal D}_ru|^2\right)\lesssim\int_0^\infty|Q_\nu^\pm(h)u|^2+\int_{3\kappa\nu}^{4\kappa\nu}\left(|u|^2+|{\cal D}_ru|^2\right).
\end{equation}
Combining (\ref{eq:3.1}) and (\ref{eq:3.2}) and taking $\lambda$ big enough, we get
\begin{equation}\label{eq:3.3}
\int_{\kappa\nu}^{2\kappa\nu}+\int_{3\kappa\nu}^{4\kappa\nu}\left(|u|^2+|{\cal D}_ru|^2\right)\lesssim\gamma (\nu+1)^2\mu_\nu\int_0^\infty(r+1)^{-2s}|{\cal D}_ru_1|^2$$
$$+\gamma^{-1}(\nu+1)^2\mu_\nu\int_0^\infty(r+1)^{2s}|Q_\nu^\pm(h)u|^2$$
$$+\varepsilon(\nu+1)^2\mu_\nu\int_0^\infty\left(|u|^2+|{\cal D}_ru|^2\right).
\end{equation}
Using (\ref{eq:3.3}) we can rewrite the above inequalities as follows:
$$\lambda\int_{5\kappa\nu/2}^\infty(r+1)^{-2s}\left(E|u|^2+|{\cal D}_ru|^2\right)$$
$$\le\lambda\int_0^\infty(r+1)^{-2s}F_1\le \int_0^\infty\psi'F_1\le\int_0^\infty\psi'\mu F_1$$
$$\lesssim\gamma\mu_\nu\int_0^\infty(r+1)^{-2s}|{\cal D}_ru_1|^2$$
$$+\gamma^{-1}\mu_\nu\int_0^\infty(r+1)^{2s}|Q_\nu^\pm(h)u|^2$$
$$+\varepsilon\mu_\nu\int_0^\infty\left(|u|^2+|{\cal D}_ru|^2\right)+\mu_\nu\int_{\kappa\nu}^{2\kappa\nu}\left(|u|^2+|{\cal D}_ru|^2\right)$$
$$\lesssim\gamma(\nu+1)^2\mu_\nu^2\int_0^\infty(r+1)^{-2s}\left(|u|^2+|{\cal D}_ru|^2\right)$$
$$+\gamma^{-1}(\nu+1)^2\mu_\nu^2\int_0^\infty(r+1)^{2s}|Q_\nu^\pm(h)u|^2$$
$$+\varepsilon(\nu+1)^2\mu_\nu^2\int_0^\infty\left(|u|^2+|{\cal D}_ru|^2\right)$$
and 
$$\int_0^{5\kappa\nu/2}(r+1)^{-2s}\left(|u|^2+|{\cal D}_ru|^2\right)\le\int_0^\infty\left(|u_2|^2+|{\cal D}_ru_2|^2\right)$$
$$\lesssim\int_0^\infty|Q_\nu^\pm(h)u|^2+\int_{3\kappa\nu}^{4\kappa\nu}
\left(|u|^2+|{\cal D}_ru|^2\right)$$
$$\lesssim\gamma(\nu+1)^2\mu_\nu\int_0^\infty(r+1)^{-2s}\left(|u|^2+|{\cal D}_ru|^2\right)$$
$$+\gamma^{-1}(\nu+1)^2\mu_\nu\int_0^\infty(r+1)^{2s}|Q_\nu^\pm(h)u|^2$$
$$+\varepsilon(\nu+1)^2\mu_\nu\int_0^\infty\left(|u|^2+|{\cal D}_ru|^2\right).$$
Hence
$$\int_0^\infty(r+1)^{-2s}\left(|u|^2+|{\cal D}_ru|^2\right)
\lesssim\gamma(\nu+1)^2\mu_\nu^2\int_0^\infty(r+1)^{-2s}\left(|u|^2+|{\cal D}_ru|^2\right)$$
$$+\gamma^{-1}(\nu+1)^2\mu_\nu^2\int_0^\infty(r+1)^{2s}|Q_\nu^\pm(h)u|^2$$
$$+\varepsilon(\nu+1)^2\mu_\nu^2\int_0^\infty\left(|u|^2+|{\cal D}_ru|^2\right).$$
We now take $\gamma$ such that $\gamma(\nu+1)^2\mu_\nu^2=\gamma_0$ with a sufficiently small constant $\gamma_0>0$, so that we can absorbe
the first term in the right-hand side of the above inequality. This leads to the estimate
$$\int_0^\infty(r+1)^{-2s}\left(|u|^2+|{\cal D}_ru|^2\right)
\lesssim(\nu+1)^4\mu_\nu^4\int_0^\infty(r+1)^{2s}|Q_\nu^\pm(h)u|^2$$
$$+\varepsilon(\nu+1)^2\mu_\nu^2\int_0^\infty\left(|u|^2+|{\cal D}_ru|^2\right),$$
which implies the desired estimate for $1/2<s\le 1$, and hence for all $s>1/2$. 
\eproof

We will next show that for large $\nu$ much better bounds for $M_\nu$ are possible. To this end, set $\tau=h^{-1/3}$, 
$\tau_1=\tau$ if
$V$ satisfies (\ref{eq:1.1}) with $\delta>2$, $\tau=h^{-\frac{1}{2\delta-1}}\epsilon^{-\frac{2-\delta}{2\delta-1}}$,
$\tau_1=\tau\epsilon^{-1}$ if $V$ satisfies (\ref{eq:1.1}) with $1<\delta\le 2$, and $\tau=h^{-1}$, $\tau_1=\tau$
if $V$ satisfies (\ref{eq:1.6}), where $\epsilon=\left(\log(h^{-1})\right)^{-1}$.

\begin{prop} There exist constants $C,c>0$ such that the estimate (\ref{eq:2.7}) holds for all 
$\nu\ge c\tau$ with $M_\nu=e^{C\tau_1/h}$.
\end{prop}

{\it Proof.} Set $\lambda=1$, $\omega(r)=(r+1)^{-2\delta+3}$ if
$V$ satisfies (\ref{eq:1.1}) with $\delta>2$, $\lambda\gg 1$, $\omega(r)=(r+1)^{-1-\epsilon}$ if
$V$ satisfies (\ref{eq:1.1}) with $1<\delta\le 2$, and $\lambda\gg 1$, $\omega(r)=(r+1)^{-1}(\log(r+2))^{-\rho}$
if $V$ satisfies (\ref{eq:1.6}). Set
$$\varphi(r)=\lambda\tau\int_0^r\omega(\sigma)d\sigma\lesssim\tau_1.$$
The parameter $\lambda$ in the second and the third cases is independent of $\nu$ and $h$ and will be fixed later on. 
Introduce the operator 
$$Q_{\nu,\varphi}^\pm(h)=e^{\varphi/h}Q_\nu^\pm(h)e^{-\varphi/h}$$
$$={\cal D}_r^2+\frac{\nu^2}{r^2}+V(r)-\varphi'^2+h\varphi''+2i\varphi'{\cal D}_r-E\pm i\varepsilon.$$
Consider the function
$$F_\varphi(r)=(E-\nu^2r^{-2}+\varphi'^2)|u(r)|^2+|{\cal D}_ru(r)|^2.$$
It is easy to see that its first derivative is given by 
$$F'_\varphi(r)=2(\nu^2r^{-3}+\varphi'\varphi'')|u|^2+4h^{-1}\varphi'|{\cal D}_ru|^2$$
$$+2h^{-1}{\rm Im}\,(V+h\varphi'')u\overline{{\cal D}_ru}-2h^{-1}{\rm Im}\,Q_{\nu,\varphi}^\pm(h)u\overline{{\cal D}_ru}\pm2\varepsilon h^{-1}{\rm Re}\,u\overline{{\cal D}_ru}$$
$$\ge \left(2\nu^2r^{-3}-2\varphi'|\varphi''|-(\varphi')^{-1}|\varphi''|^2\right)|u|^2-h^{-1}(\varphi')^{-1}V^2|u|^2-\Phi(r),$$
where
$$\Phi=\gamma(r+1)^{-2s}\left|{\cal D}_ru\right|^2+\gamma^{-1}h^{-2}(r+1)^{2s}\left|Q_{\nu,\varphi}^\pm(h)u\right|^2$$
$$+\varepsilon h^{-1}\left(|u|^2+|{\cal D}_ru|^2\right)$$
and $\gamma>0$ and $s>1/2$ are arbitrary. 
 It is easy to check that
$$2\varphi'|\varphi''|+(\varphi')^{-1}|\varphi''|^2\le c_0^2\tau^2(r+1)^{-3}$$
with some constant $c_0>0$. Therefore, for $\nu\ge c_0\tau$, we obtain
\begin{equation}\label{eq:3.4}
F'_\varphi(r)\ge \nu^2r^{-3}|u|^2-(\lambda h\tau)^{-1}\omega^{-1}V^2|u|^2-\Phi.
\end{equation}
Let us see that (\ref{eq:3.4}) implies the inequality
\begin{equation}\label{eq:3.5}
F'_\varphi(r)\ge \frac{1}{2}\nu^2r^{-3}|u|^2-K\lambda^{-1}\omega|u|^2-\Phi,
\end{equation}
provided $\nu\ge c\tau$ with some constant $c>c_0$, where $K=0$ if
$V$ satisfies (\ref{eq:1.1}) with $\delta>2$, $K>0$ is a constant if
$V$ satisfies (\ref{eq:1.6}), and $K=\widetilde K\epsilon$, $\widetilde K>0$ is a constant if
$V$ satisfies (\ref{eq:1.1}) with $1<\delta\le 2$. Indeed, in the first case this follows from
the inequality $(h\tau)^{-1}\omega^{-1}V^2\lesssim \tau^2(r+1)^{-3}$, while in the case when 
$V$ satisfies (\ref{eq:1.6}) it follows from $(h\tau)^{-1}\omega^{-1}V^2\lesssim \omega$.
Let now $V$ satisfy (\ref{eq:1.1}) with $1<\delta\le 2$. Then we have
$$\omega^{-1}V^2\lesssim (r+1)^{-2\delta+1+\epsilon}\le a^{4-2\delta+\epsilon}(r+1)^{-3}+a^{2-2\delta+2\epsilon}(r+1)^{-1-\epsilon}$$
for every $a>1$. Take $a$ such that
$$(h\tau)^{-1}a^{4-2\delta}=\tau^2,\quad (h\tau)^{-1}a^{2-2\delta}=\epsilon.$$
In view of the choice of $\tau$ we find that these equations are satisfied with 
$a=h^{-\frac{1}{2\delta-1}}\epsilon^{-\frac{3}{2(2\delta-1)}}$.
Thus we get the inequality
$$(h\tau)^{-1}\omega^{-1}V^2\lesssim \tau^2(r+1)^{-3}+\epsilon(r+1)^{-1-\epsilon},$$
which clearly implies (\ref{eq:3.5}) in this case. 

Integrating (\ref{eq:3.5}) from $0$ to $\infty$ and using that $F_\varphi(0)=|{\cal D}_ru(0)|^2\ge 0$, we get
\begin{equation}\label{eq:3.6}
\frac{1}{2}\int_0^\infty\nu^2r^{-3}|u|^2\le K\lambda^{-1}\int_0^\infty\omega|u|^2+\int_0^\infty\Phi.
\end{equation}
On the other hand, integrating (\ref{eq:3.5}) from $r$ to $\infty$ yields
\begin{equation}\label{eq:3.7}
F_\varphi(r)=-\int_r^\infty F'_\varphi(\sigma)d\sigma\le 
K\lambda^{-1}\int_0^\infty\omega|u|^2+\int_0^\infty\Phi.
\end{equation}
Set $\widetilde\omega=(r+1)^{-2s}$ if
$V$ satisfies (\ref{eq:1.1}) with $\delta>2$, and $\widetilde\omega=\omega$ otherwise. Clearly, 
$\|\widetilde\omega\|_{L^1}\lesssim \epsilon^{-1}$ if $V$ satisfies (\ref{eq:1.1}) with $1<\delta\le 2$
and $\|\widetilde\omega\|_{L^1}\lesssim 1$ in the other two cases. Therefore, multiplying (\ref{eq:3.7}) by $\widetilde\omega$ and integrating
 from $0$ to $\infty$ lead to
\begin{equation}\label{eq:3.8}
\int_0^\infty\widetilde\omega F_\varphi\lesssim  
\ell_1\lambda^{-1}\int_0^\infty\omega|u|^2+(1+\ell_2\epsilon^{-1})\int_0^\infty\Phi,
\end{equation}
where $\ell_1=\ell_2 =0$ if $V$ satisfies (\ref{eq:1.1}) with $\delta>2$, $\ell_1=\ell_2=1$ if $V$ satisfies (\ref{eq:1.1}) with $1<\delta\le 2$,
and $\ell_1=1,\ell_2=0$ if $V$ satisfies (\ref{eq:1.6}). By (\ref{eq:3.6}) and (\ref{eq:3.8}), using that $\widetilde\omega\le Cr^{-1}$, we obtain 
$$\int_0^\infty\widetilde\omega\left(E|u|^2+|{\cal D}_ru|^2\right)\le 
\int_0^\infty\widetilde\omega F_\varphi+C\int_0^\infty \nu^2r^{-3}|u|^2$$
$$\lesssim  
\ell_1\lambda^{-1}\int_0^\infty\widetilde\omega|u|^2+(1+\ell_2\epsilon^{-1})\int_0^\infty\Phi.$$
When $\ell_1=1$ we take $\lambda$ large enough in order to absorbe the first term in the right-hand side of the above inequality. Thus we get
\begin{equation}\label{eq:3.9}
\int_0^\infty\widetilde\omega\left(|u|^2+|{\cal D}_ru|^2\right)
\lesssim (1+\ell_2\epsilon^{-1})\int_0^\infty\Phi.
\end{equation}
Since $\widetilde\omega\ge (r+1)^{-2s}$, we deduce from (\ref{eq:3.9})
$$\int_0^\infty(r+1)^{-2s}\left(|u|^2+|{\cal D}_ru|^2\right)
\lesssim \gamma(1+\ell_2\epsilon^{-1})\int_0^\infty(r+1)^{-2s}\left|{\cal D}_ru\right|^2$$
$$+\gamma^{-1}h^{-2}(1+\ell_2\epsilon^{-1})\int_0^\infty(r+1)^{2s}\left|Q_{\nu,\varphi}^\pm(h)u\right|^2$$
$$+\varepsilon h^{-1}(1+\ell_2\epsilon^{-1})\int_0^\infty\left(|u|^2+|{\cal D}_ru|^2\right).$$
We now take $\gamma$ such that $\gamma(1+\ell_2\epsilon^{-1})=\gamma_0$, where $\gamma_0$ is a sufficiently small constant.
Thus we can absorbe the first term in the right-hand side of the above inequality to obtain
\begin{equation}\label{eq:3.10}
\int_0^\infty(r+1)^{-2s}|u|^2
\lesssim h^{-2}(1+\ell_2\epsilon^{-1})^2\int_0^\infty(r+1)^{2s}\left|Q_{\nu,\varphi}^\pm(h)u\right|^2$$ 
$$+\varepsilon h^{-1}(1+\ell_2\epsilon^{-1})\int_0^\infty\left(|u|^2+|{\cal D}_ru|^2\right).
\end{equation}
We apply (\ref{eq:3.10}) with $u$ replaced by $e^{\varphi/h}u$. Thus we get the Carleman estimate
\begin{equation}\label{eq:3.11}
\int_0^\infty(r+1)^{-2s}e^{2\varphi/h}|u|^2
\lesssim h^{-2}(1+\ell_2\epsilon^{-1})^2\int_0^\infty(r+1)^{2s}e^{2\varphi/h}\left|Q_{\nu}^\pm(h)u\right|^2$$ 
$$+\varepsilon h^{-1}(1+\ell_2\epsilon^{-1})\tau^2\int_0^\infty e^{2\varphi/h}\left(|u|^2+|{\cal D}_ru|^2\right).
\end{equation}
Since $1\le e^{2\varphi/h}\le e^{\widetilde C\tau_1/h}$ with some constant $\widetilde C>0$, (\ref{eq:3.11}) implies
\begin{equation}\label{eq:3.12}
\int_0^\infty(r+1)^{-2s}|u|^2
\lesssim h^{-2}(1+\ell_2\epsilon^{-1})^2e^{\widetilde C\tau_1/h}\int_0^\infty(r+1)^{2s}\left|Q_{\nu}^\pm(h)u\right|^2$$ 
$$+\varepsilon h^{-1}(1+\ell_2\epsilon^{-1})\tau^2e^{\widetilde C\tau_1/h}\int_0^\infty\left(|u|^2+|{\cal D}_ru|^2\right),
\end{equation}
which gives the desired bound for $M_\nu$.
\eproof

It follows from Propositions 3.1 and 3.2 together with Lemma 2.2 that $M\le e^{C\tau_1/h}$ with some constant $C>0$, 
which, in view of Lemma 2.1, implies Theorem 1.1.

\section{Bounding $M$ for H\"older potentials} 

Let $\rho\in C_0^\infty([0,1])$, $\rho\ge 0$, be a real-valued function independent of $\nu$ and $h$ such that $\int_0^\infty\rho(\sigma)d\sigma=1$.
If $V\in C^\alpha_\beta(\mathbb{R}^+)$, we can approximate it by the function
$$V_\theta(r)=\theta^{-1}\int_0^\infty\rho((r-r')/\theta)V(r')dr'=\int_0^\infty\rho(\sigma)V(r+\theta\sigma)d\sigma$$
where $0<\theta\ll 1$ will be chosen later on. We have
\begin{equation}\label{eq:4.1}
|V(r)-V_\theta(r)|\le\int_0^\infty\rho(\sigma)|V(r+\theta\sigma)-V(r)|d\sigma$$
$$\lesssim\theta^\alpha(r+1)^{-\beta}\int_0^\infty\sigma^\alpha\rho(\sigma)d\sigma\lesssim\theta^\alpha(r+1)^{-\beta}.
\end{equation}
If in addition $V$ satisfies (\ref{eq:1.9}), we obtain from (\ref{eq:4.1})
\begin{equation}\label{eq:4.2}
V_\theta(r)\le p(r)+{\cal O}((r+1)^{-\beta}).
\end{equation}
Since $p$ is a decreasing function tending to zero, it follows from (\ref{eq:4.2})
that given any $\gamma, N>0$ there is a constant $C_{\gamma,N}>0$ such that
\begin{equation}\label{eq:4.3}
V_\theta(r)\le C_{\gamma,N}(r+1)^{-N}+\gamma.
\end{equation}
Clearly, $V_\theta\in C^1$ and its first derivative $V'_\theta$ is given by
$$V'_\theta(r)=\theta^{-2}\int_0^\infty\rho'((r-r')/\theta)V(r')dr'$$
$$=\theta^{-1}\int_0^\infty\rho'(\sigma)V(r+\theta\sigma)d\sigma=\theta^{-1}\int_0^\infty\rho'(\sigma)(V(r+\theta\sigma)-V(r))d\sigma$$
where we have used that $\int_0^\infty\rho'(\sigma)d\sigma=0$. Hence
\begin{equation}\label{eq:4.4}
|V'_\theta(r)|\lesssim\theta^{-1+\alpha}(r+1)^{-\beta}\int_0^\infty\sigma^\alpha|\rho'(\sigma)|d\sigma
\lesssim\theta^{-1+\alpha}(r+1)^{-\beta}.
\end{equation}
Using the above inequalities we will prove the following

\begin{prop} Let $V\in L^\infty$ satisfy (\ref{eq:1.9}). Suppose in addition that
$V\in C^\alpha_\beta(\mathbb{R}^+)$ with $\beta>1$. 
Then the estimate (\ref{eq:2.7}) holds for all 
$\nu$ with $M_\nu=e^{C(\nu+1)/h}$, where $C>0$ is a constant independent of $\nu$ and $h$.
\end{prop}

{\it Proof.} We will modify the proof of Proposition 3.1 to avoid using that $V\in L^1$. Instead, we will use that
$V-V_\theta$ and $V'_\theta$ belong to $L^1$. 
We will apply the above inequalities 
with $\theta$ independent of $h$ and $\nu$. Set
$\varphi(r)=\lambda\left(1-(r+1)^{-1}\right)$, where $\lambda\gg 1$ is independent of $h$ and $\nu$. 
Clearly, $\varphi'(r)=\lambda(r+1)^{-2}$. 
Using (\ref{eq:4.3}) with $\gamma=E/2$ and $N=4$, we obtain
\begin{equation}\label{eq:4.5}
E+\varphi'^2-V_\theta\ge \frac{E}{2}+(\lambda^2-C_{\gamma,4})(r+1)^{-4}\ge \frac{E}{2}
\end{equation}
provided $\lambda$ is large enough. 
We set 
$$F(r)=(E+\varphi'^2-V_\theta)|u(r)|^2+|{\cal D}_ru(r)|^2$$
when $\nu=0$, and 
$$F_1(r)=(E+\varphi'^2-V_\theta)|u_1(r)|^2+|{\cal D}_ru_1(r)|^2$$
when $\nu>0$. By (\ref{eq:4.5}),
\begin{equation}\label{eq:4.6}
\frac{E}{2}|u(r)|^2+|{\cal D}_ru(r)|^2\le F(r).
\end{equation}
Clearly, (\ref{eq:4.6}) holds with $u$ and $F$ replaced by $u_1$ and $F_1$. 
As in the previous section, the first derivative of $F$ is given by
$$F'(r)=(2\varphi'\varphi''-V'_\theta)|u|^2+4h^{-1}\varphi'|{\cal D}_ru|^2$$
$$+2h^{-1}{\rm Im}\,(V-V_\theta+h\varphi'')u\overline{{\cal D}_ru}-2h^{-1}{\rm Im}\,Q_{0,\varphi}^\pm(h)u\overline{{\cal D}_ru}\pm2\varepsilon h^{-1}{\rm Re}\,u\overline{{\cal D}_ru}.$$
Hence
$$-F'(r)\le h^{-1}W_0\left(|u|^2+|{\cal D}_ru|^2\right)+h^{-1}\gamma(r+1)^{-2s}|{\cal D}_ru|^2$$
$$+h^{-1}\gamma^{-1}(r+1)^{2s}|Q_{0,\varphi}^\pm(h)u|^2
+\varepsilon h^{-1}\left(|u|^2+|{\cal D}_ru|^2\right)$$
for any $\gamma>0$, where 
$$W_0=|V-V_\theta|+|V'_\theta|+|\varphi''|+2\varphi'|\varphi''|\in L^1.$$
Similarly
$$-F'_1(r)\le h^{-1}W_\nu\left(|u_1|^2+|{\cal D}_ru_1|^2\right)+h^{-1}\gamma(r+1)^{-2s}|{\cal D}_ru_1|^2$$
$$+h^{-1}\gamma^{-1}(r+1)^{2s}|Q_{\nu,\varphi}^\pm(h)u|^2
+\varepsilon{\cal O}\left(h^{-1}\right)\left(|u|^2+|{\cal D}_ru|^2\right)$$
$$+{\cal O}\left(h\nu^{-2}\right)(\phi'_1+|\phi''_1|)\left(|u|^2+|{\cal D}_ru|^2\right)$$
for any $\gamma>0$, where
$$W_\nu=\nu^2r^{-2}\phi(r/\kappa\nu)+W_0\in L^1.$$
Clearly, $\|W_\nu\|_{L^1}={\cal O}(\nu+1)$. 
  Now, arguing in the same way as in the proof of Proposition 3.1 and using
 the above inequalities, we obtain the estimate
 $$\int_0^\infty(r+1)^{-2s}|u|^2\le e^{C(\nu+1)/h}\int_0^\infty(r+1)^{2s}|Q_{\nu,\varphi}^\pm(h)u|^2
+\varepsilon e^{C(\nu+1)/h}\int_0^\infty\left(|u|^2+|{\cal D}_ru|^2\right)$$
with some constant $C>0$. Applying this inequality with $u$ replaced by $e^{\varphi/h}u$, we get
$$\int_0^\infty(r+1)^{-2s}e^{2\varphi/h}|u|^2$$
$$\le e^{C(\nu+1)/h}\int_0^\infty(r+1)^{2s}e^{2\varphi/h}|Q_\nu^\pm(h)u|^2
+\varepsilon e^{C(\nu+1)/h}\int_0^\infty e^{2\varphi/h}\left(|u|^2+|{\cal D}_ru|^2\right)$$
which implies
$$\int_0^\infty(r+1)^{-2s}|u|^2\le e^{C(\nu+1)/h}\int_0^\infty(r+1)^{2s}|Q_\nu^\pm(h)u|^2
+\varepsilon e^{C(\nu+1)/h}\int_0^\infty\left(|u|^2+|{\cal D}_ru|^2\right)$$
with a new constant $C>0$, which is the desired estimate.
\eproof

Let $V\in C_\beta^\alpha$ with $0<\alpha<1$ and $1<\beta\le 3$. Set $\tau=h^{-k_0}\epsilon^{-q_0}$, $\tau_1=h^{-k_0}\epsilon^{-q},$
where
$$k_0=\frac{1-\alpha}{\alpha+3},\quad q=q_0=0,$$ 
if $\beta=3$,
$$k_0=\frac{1-\alpha}{2\alpha\beta-5\alpha+3},\quad q=q_0=\frac{\alpha(3-\beta)}{2\alpha\beta-5\alpha+3},$$
  if $2<\beta<3$, 
$$k_0=\frac{1-\alpha}{2\beta-\alpha-1},\quad q_0=
\frac{\alpha-\beta+2}{2\beta-\alpha-1},\quad q=q_0+1,$$ 
  if $1<\beta\le 2$. Clearly, to get Theorem 1.2 it suffices to prove the following

\begin{prop} There exist constants $C,c>0$ such that the estimate (\ref{eq:2.7}) holds for all 
$\nu\ge c\tau$ with $M_\nu=e^{C\tau_1/h}$.
\end{prop}

{\it Proof.} Set $\omega(r)=(r+1)^{-2\beta+3}$ if $2<\beta\le 3$, $\omega(r)=(r+1)^{-1-\epsilon}$ if $1<\beta\le 2$ and
$$\varphi(r)=\lambda\tau\int_0^r\omega(\sigma)d\sigma\lesssim\tau_1,$$
where $\lambda\gg 1$ is independent of $\nu$ and $h$. Consider the function
$$F_\varphi(r)=(E-\nu^2r^{-2}+\varphi'^2-V_\theta)|u(r)|^2+|{\cal D}_ru(r)|^2.$$
Clearly, we can still arrange the inequality (\ref{eq:4.5}). Therefore, we have
\begin{equation}\label{eq:4.7}
\frac{E}{2}|u(r)|^2+|{\cal D}_ru(r)|^2\le F_\varphi(r)+\nu^2r^{-2}|u(r)|^2.
\end{equation}
In the same way as in the proof of Proposition 3.2 we can obtain the following analog of the inequality (\ref{eq:3.4})
\begin{equation}\label{eq:4.8}
F'_\varphi(r)\ge \nu^2r^{-3}|u|^2-W|u|^2-\Phi,
\end{equation}
for $\nu\ge c_0\tau$, where
$$W=|V'_\theta|+(\lambda h\tau)^{-1}\omega^{-1}|V-V_\theta|^2.$$
By (\ref{eq:4.1}) and (\ref{eq:4.4}), we have
\begin{equation}\label{eq:4.9}
W\lesssim \theta^{-1+\alpha}(r+1)^{-\beta}+(\lambda h\tau)^{-1}\theta^{2\alpha}\omega^{-1}(r+1)^{-2\beta}.
\end{equation}
Set $\ell=0$, $\widetilde\omega=(r+1)^{-2s}$ if $\beta=3$ and $\ell=1$, $\widetilde\omega=(r+1)^{-1-\epsilon}$ if $1<\beta<3$.
We will show that (\ref{eq:4.8}) and (\ref{eq:4.9}) imply
\begin{equation}\label{eq:4.10}
F'_\varphi(r)\ge \frac{1}{2}\nu^2r^{-3}|u|^2-{\cal O}(\lambda^{-1})\ell\epsilon\widetilde\omega|u|^2-\Phi
\end{equation}
for $\nu\ge c\tau$ with some constant $c\gg 1$. If $\beta=3$ we require that $\tau$ and $\theta$ satisfy the relations
$$\tau^2=\theta^{-1+\alpha}=(h\tau)^{-1}\theta^{2\alpha},$$
which provides the desired value of $\tau$. Since in this case $\omega^{-1}(r+1)^{-2\beta}=(r+1)^{-3}$, we get 
(\ref{eq:4.10}) from (\ref{eq:4.8}) and (\ref{eq:4.9}). Let $2<\beta<3$. Then we have the inequality
\begin{equation}\label{eq:4.11}
(r+1)^{-\beta}\le (bb_0)^{3-\beta}(r+1)^{-3}+(bb_0)^{1+\epsilon-\beta}(r+1)^{-1-\epsilon}
\end{equation}
for every $b,b_0>1$, provided $\epsilon\ll \beta-1$. We choose $b_0$ such that $b_0^{1-\beta}=\lambda^{-1}$.
We also let 
$$\tau^2=b^{3-\beta}\theta^{-1+\alpha}=(h\tau)^{-1}\theta^{2\alpha},\quad \theta^{-1+\alpha}b^{-\beta+1}=\epsilon.$$
We are looking for solutions of these equations of the form $\tau=h^{-k_0}\epsilon^{-q_0}$,
$\theta^{-1}=h^{-k_1}\epsilon^{-q_1}$, $b=h^{-k_2}\epsilon^{-q_2}$. Thus the above
equations take the form
\begin{equation}\label{eq:4.12}
\left\{
\begin{array}{l}
2k_0=(3-\beta)k_2+(1-\alpha)k_1=1-k_0-2\alpha k_1,\\
(1-\alpha)k_1-(\beta-1)k_2=0,\\
2q_0=(3-\beta)q_2+(1-\alpha)q_1=-q_0-2\alpha q_1,\\
(1-\alpha)q_1-(\beta-1)q_2=-1.
\end{array}
\right.
\end{equation}
Solving this linear system we find the desired values of $k_0$ and $q_0$. 
With this choice, by (\ref{eq:4.9}) we get
\begin{equation}\label{eq:4.13}
W\lesssim \tau^2(r+1)^{-3}+\lambda^{-1}\epsilon(r+1)^{-1-\epsilon},
\end{equation}
which together with (\ref{eq:4.8}) imply (\ref{eq:4.10}). 

Let $1<\beta<2$. Then (\ref{eq:4.11}) still holds. We also have the inequality
\begin{equation}\label{eq:4.14}
\omega^{-1}(r+1)^{-2\beta}=(r+1)^{-2\beta+1+\epsilon}\le a^{4-2\beta+\epsilon}(r+1)^{-3}+a^{2-2\beta+2\epsilon}(r+1)^{-1-\epsilon}
\end{equation}
for every $a>1$. We let the parameters $a$ and $b$ satisfy the relations
$$\tau^2=b^{3-\beta}\theta^{-1+\alpha}=(h\tau)^{-1}\theta^{2\alpha}a^{4-2\beta},\quad 
\theta^{-1+\alpha}b^{-\beta+1}=(h\tau)^{-1}\theta^{2\alpha}a^{2-2\beta}=\epsilon.$$
As above, we are looking for solutions of these equations of the form $\tau=h^{-k_0}\epsilon^{-q_0}$,
$\theta^{-1}=h^{-k_1}\epsilon^{-q_1}$, $b=h^{-k_2}\epsilon^{-q_2}$ and $a=h^{-k_3}\epsilon^{-q_3}$. Thus the above
equations take the form
\begin{equation}\label{eq:4.15}
\left\{
\begin{array}{l}
2k_0=(3-\beta)k_2+(1-\alpha)k_1=1-k_0-2\alpha k_1+2(2-\beta)k_3,\\
(1-\alpha)k_1-(\beta-1)k_2=1-k_0-2\alpha k_1-2(\beta-1)k_3=0,\\
2q_0=(3-\beta)q_2+(1-\alpha)q_1=-q_0-2\alpha q_1+2(2-\beta)q_3,\\
(1-\alpha)q_1-(\beta-1)q_2=-q_0-2\alpha q_1-2(\beta-1)q_3=-1.
\end{array}
\right.
\end{equation}
Again, solving this linear system we find the desired values of $k_0$ and $q_0$. Thus we conclude that
(\ref{eq:4.13}) still holds in this case, and hence (\ref{eq:4.10}) follows.
On the other hand, (\ref{eq:4.7}) and (\ref{eq:4.10}) imply the estimate (\ref{eq:2.7})
in the same way as in the proof of Proposition 3.2.
\eproof


\begin{thebibliography}
\frenchspacing \baselineskip=12 pt plus 1pt minus 1pt

\bibitem{kn:B1} {\sc N. Burq}, {\em D\'ecroissance de l'\'energie locale de l'\'equation des ondes pour le probl\`eme
ext\'erieur et absence de r\'esonance au voisinage du r\'eel}, Acta Math. {\bf 180} (1998), 1-29.

\bibitem{kn:B2} {\sc N. Burq}, {\em Lower bounds for shape resonances widths of long-range Schr\"odinger operators},
Amer. J. Math. {\bf 124} (2002), 677-735.

\bibitem{kn:CV} {\sc F. Cardoso and G. Vodev}, {\em Uniform estimates of the resolvent of the Laplace-Beltrami operator
on infinite volume Riemannian manifolds}, Ann. Henri Poincar\'e {\bf 4} (2002), 673-691.

\bibitem{kn:GS} {\sc J. Galkowski and J. Shapiro}, {\em Semiclassical resolvent bounds for weakly decaying potentials}, preprint 2020.

\bibitem{kn:D} {\sc K. Datchev}, {\em Quantative limiting absorption principle in the semiclassical limit}, Geom. Funct. Anal.
{\bf 24} (2014), 740-747.

\bibitem{kn:DDZ} {\sc K. Datchev, S. Dyatlov and M. Zworski}, {\em Resonances and lower resolvent bounds}, J. Spectral Theory
{\bf 5} (2015), 599-615.

\bibitem{kn:DS} {\sc K. Datchev and J. Shapiro}, {\em Semiclassical estimates for scattering on the real line}, Commun. Math.
Phys. {\bf 376} (2020), 2301-2308.

\bibitem{kn:KV} {\sc F. Klopp and M. Vogel}, {\em Semiclassical resolvent estimates for bounded potentials}, 
Pure Appl. Analysis {\bf 1} (2019), 1-25.

\bibitem{kn:S1} {\sc J. Shapiro}, {\em Local energy decay for Lipschitz wavespeeds}, Commun. Partial Diff. Equations {\bf 43} (2018), 839-858.

\bibitem{kn:S2} {\sc J. Shapiro}, {\em Semiclassical resolvent bounds in dimension two}, Proc. Amer. Math. Soc. {\bf 147} (2019), 1999-2008.

\bibitem{kn:S3} {\sc J. Shapiro}, {\em Semiclassical resolvent bound for compactly supported $L^\infty$ potentials},
J. Spectral Theory {\bf 10} (2020), 651-672.

\bibitem{kn:V1} {\sc G. Vodev}, {\em Semiclassical resolvent estimates for short-range $L^\infty$ potentials}, 
Pure Appl. Analysis {\bf 1} (2019), 207-214.

\bibitem{kn:V2} {\sc G. Vodev}, {\em Semiclassical resolvent estimates for short-range $L^\infty$ potentials. II}, Asymptotic Analysis, to appear.

\bibitem{kn:V3} {\sc G. Vodev}, {\em Semiclassical resolvent estimates for $L^\infty$ potentials
on Riemannian manifolds}, Ann. Henri Poincar\'e {\bf 21} (2020), 437-459.

\bibitem{kn:V4} {\sc G. Vodev}, {\em Semiclassical resolvent estimates for H\"older potentials}, preprint 2020.

\end{thebibliography}
\end{document}